\documentclass[12pt]{article}
\usepackage{amsmath,amssymb,amsfonts}

\newcommand{\var}{\varepsilon}
\newcommand{\rg}{\rightarrow}
\newcommand{\R}{\mathbb{{R}}}
\newcommand{\para}{\paragraph}

\begin{document}
\begin{center}
{\bf Multiple Equilibria for an SIRS Epidemiological System}

\vspace*{1cm}
M.R. Razvan\\
{\it\small Institute for Studies in Theoretical Physics and Mathematics\\
P.O. Box 19395-5746, Tehran, Iran\\
e-mail: razvan@karun.ipm.ac.ir}

\begin{abstract}
An SIRS type model of disease transmission in an open environment
is discussed. We use the Poincar\'{e} index together with a
perturbation method to show that the endemic proportions need not
be unique.
\end{abstract}
\end{center}

\section{Introduction}
The social mixing structure of a population or a group of
interacting populations play a crucial role in the dynamics of a
disease transmission \cite{CC}.  The most well-known examples of
epidemics which are spread by means of the interaction between
different populations are those related to venereal diseases.
These diseases are transmitted by sexual contacts between two or
more different populations. An epidemiological model which
considers these interactions is called multigroup model. A large
number of multigroup models have been described in \cite{C}.
These models usually lead to some high-dimensional systems of
differential equations with a probably complicated dynamics
\cite{HC}. In order to avoid these complications, we consider a
single group and we assume that there is a disease transmission
from the outside. This assumption is natural, since we are living
in an open environment with a limited control on the outside. The
effect of this assumption has been studied by using numerical and
statistical methods. For example in \cite{AM}, it has been
claimed that most of Swedish HIV-patients with heterosexual
transmission had been in contact with a partner of foreign
extraction. The effect of the immigrant subpopulation has been
considered by many others too \cite{A,CK,DD}. A statistical
approach to this concept in \cite{CS} showed that the
subpopulation with a high prevalence of HIV-infected individuals
is the immigrant one.

In this paper, we consider a simple SIRS type model of disease
transmission in an open environment. In an SIRS model, we divide
the population into three classes consisting of  Susceptible,
Infective and Removed individuals. By an open environment we mean
that our population has contacts with other populations and there
is a disease transmission from the outside \cite{HM}. We assume
that a proportion of susceptible individuals are infected in this
way. We also assume that there is a special control on the people
who are frequently in contact with foreigners and a proportion of
those who are infected in this way are known and removed. In our
model, the incidence function is of proportionate mixing type
introduced by Nold \cite{N}. The demographic assumptions are also
very simple. Natural births and deaths are assumed to be
proportional to the class numbers with all newborns susceptibles
and the excess deaths due to the disease among infectives and
removeds are proportional too. We could have considered more
complicated demographic or vertical transmission assumptions
\cite{BC,BH}, but these parameters have no mathematical
significance and one can easily conclude that our main results
would still be valid. We want to avoid these complications for two
reasons. The first one is that even in this simple case, multiple
endemic equilibria may occur. Indeed we shall show that for some
suitable values of parameters, the proportions system admits two
sinks and one saddle point in the feasibility region. The
existence of multiple equilibria or limit cycles has been shown
for more complicated systems such as multigroup models \cite{HC}
or a single group with nonlinear incidence functions
\cite{Dv,HL,Hv}. The second reason is related to the technique
used here to determine the number of endemic equilibria. We state
a simple model in order that we can exhibit our technique more
clearly. The technique used here is based on a careful choice of
Jordan curves and counting the number of rest points inside them.
This technique has no hard analysis and can be easily applied to
other similar systems \cite{R,RK}. The reader can verify that our
results hold for similar SIR and SIRI systems as well.

We first in the next section, state the model and some results
concerning the non-existence of certain types of solutions. In
Section 3, we use Poincar\'{e} index to  obtain some partial
results in a special case. Then in Section 4, we use these
results together with a perturbation method to examine the
general case. Our analysis relies on an index lemma concerning
the Poincar\'e index of a class of Jordan curves. We provide the
proof of this result at the end of this paper.

\section{The Model}

We consider a model of disease transmission in as nonconstant
population of size $N$ divided into three classes: susceptibles,
infectives and removeds, the number of each class is given by
$S,\ I,\ R$, respectively.  We set $N=S+I+R$ and use the
following parameters which are assumed to be positive unless
otherwise specified:
\begin{itemize}
\item[$b$=] per capita birth rate,
\item[$d$=] per capita disease free death rate,
\item[$\var_1$=] excess per capita death rate of infectives,
\item[$\var_2$=] excess per capita death rate of removeds,
\item[$\lambda$=] effective per capita contact rate of infectives,
\item[$\alpha$=] per capita removal rate of infectives,
\item[$\gamma$=] per capita recovery rate of removeds.
\end{itemize}
Here as mentioned before, we assume that there is a disease
transmission from the outside to our population. Let $\beta >0$
be the per capita transfer rate of the disease into this
population. Furthermore the susceptible individuals which are
infected in this way, enter to the classes $I$ or $R$ of
proportions $\beta_1$ and $\beta_2$ where
$\beta_1+\beta_2=\beta$. The above hypotheses lead to the
following system of differential equations in $\R^3_+$, where ``
$'$ " denotes the derivative with respect to $t$, the time.
$$\left\{
\begin{array}{lrrr}
S'=bN-(d+\beta )S+\gamma R-\dfrac{\lambda I S}{N},&&&(2-1)
\\
I'=\beta_1S-(d+\var_1+\alpha )I+\dfrac{\lambda I S}{N},
&&&(2-2)\\
R'=\beta_2S-(d+\var_2+\gamma )R+\alpha I,&&&(2-3)
\end{array}\right.$$
where $\dfrac{\lambda I}{N}$ is of the proportionate mixing type
introduced by Nold \cite{N}. The total population equation is
obtained by adding the above three equations:
$$N'=(b-d)N-\var_1 I-\var_2R. \eqno{(2-4)} $$
Now if we set $s=\dfrac{S}{N}, i=\dfrac{I}{N}$ and
$r=\dfrac{R}{N}$, the equations (2-1)-(2-3) yield
$$\left\{
\begin{array}{llrrr}
s'=&b-(b+\beta )s+\gamma r+(\var_1-\lambda )i s+
\var_2rs, &&&(2-1)'\\
i'=&\beta_1 s-(b+\var_1+\alpha )i+\lambda i s +\var_1 i^2+
\var_2 i r, &&&(2-2)'\\
r'=&\beta_2 s-(b+\var_2+\gamma )r+\alpha i +\var_1 i r +\var_2
r^2. &&&(2-3)'
\end{array}\right.$$
We will determine the asymptotic behaviour of the solutions of
this system. Prior to this we have the following  concepts of
ODE's related to our system.

Given an autonomous system of ordinary differential equations in
$\R^n$

$$\frac{dx}{dt}=f(x) \eqno{(2-5)}$$
We will denote by $x.t$ the value of the solution of this system
at time $t$ that is $x$ initially. For $V\subset\R^n,
J\subset\R$, we let $V.J=\{ x.t: x\in V, t\in J\}$. The set $V$
is called invariant if $V.\R=V$ and it is called positively
invariant if $V.\R_+=V$. For $Y\subset\R^n$ the $\omega$-limit
set of $Y$ is defined to be the maximal invariant set in the
closure of $Y.[0,\infty )$. A closed curve connecting several
rest points whose segments between successive rest points are
heteroclinic orbits is called phase polygon. By a sink we mean a
rest point at which all the eigenvalues of the linearized system
have negative real parts. A rest point is called a source if
these eigenvalues have positive real parts and it is called a
saddle point if some of these eigenvalues have positive real
parts and the others have negative real parts. A rest point is
called nondegenerate if all of these eigenvalues are nonzero and
it is called hyperbolic if all of its eigenvalues have nonzero
real parts.

Now we continue the analysis of the system $(2-1)'-(2-3)'$.
 If we set $\sum =s +i+r$ then $\sum^{\prime}
=(1-\sum )(b-\var_1 i-\var_2 r)$. Therefore the plane $\sum=1$ is
invariant. We consider the feasibility region
$$D=\{ (s, i, r)|s+i+r=1, s\geq 0, i\geq 0, r\geq 0\},$$
which is a triangle and on its sides we have
$$
s=0\Rightarrow s'=b+\gamma r \ \ ,\ \ i=0\Rightarrow i'=\beta_1 s
\ \ ,\ \  r=0\Rightarrow r'=\beta_2 s+\gamma i.$$ Since all
parameters are positive,
 $D$ is positively invariant and any solutions of the system
$(2-1)'-(2-3)'$ with initial point in $\partial D$ immediately
enters $\overset{\circ}{D}$, where $\partial D$ and
$\overset{o}{D}$ are the boundary and the interior of $D$,
respectively. From now on, we examine the dynamics of this system
in the feasibility region $D$.

Using the relation $s+i+r=1$, we see that the system is essentially two
dimensional. Thus we
 can eliminate one of the variables to arrive at the following quadratic
planner system: $$\left\{
\begin{array}{lr}
s'=b+\gamma +(\var_2-b-\beta -\gamma )s-\gamma i-
\var_2 s^2+(\var_1-\var_2-\lambda )i s,&(2-6)\\
i'=\beta_1 s+(\var_2-\var_1-\alpha -b)i+(\lambda -\var_2) i
s+(\var_1-\var_2)i^2. &(2-7)
\end{array}
\right. $$ Notice that this planner system has at most four rest
points. Moreover the dynamics of the system $(2-1)'-(2-3)'$ on $D$
is equivalent to the dynamics of this system in the positively
invariant region $D_1=\{ (s, i)|s\geq 0, i\geq 0, s+i\leq 1 \}$.

The following theorem is a special case of the results of
\cite{Bv2} concerning the non-existence of certain types of
solutions.

\para{Theorem 2.1.} Let in (2-5), $f$
 be a smooth vector field in $\R^3$. Let $\Gamma
(t)$ be a closed piece-wise smooth curve, which is the boundary of
an orientable smooth surface $S\subset\R^3$. Suppose
$g:U\rightarrow\R^3$ is defined and smooth in a neighborhood $U$
of $S$ with $g(\Gamma (t)).f(\Gamma (t))\geq 0$ and
$\mbox{curl}(g).n<0$, where $n$ is the unit normal to $S$. Then
$\Gamma$ is not a finite union of the orbits of the system (2.5).

\para{Proposition 2.2.} The system $(2-1)'-(2-3)'$ has no periodic orbits,
homoclinic orbits, or phase polygons in $\overset{o}{D}$.\\
{\bf Proof.} In order to apply the above theorem, we define
 $g=g_1+g_2+g_3$ where
\begin{align*}
g_1(i, r)=\left [0, -\dfrac{f_3(i, r)}{ir}, \dfrac{f_2(i, r)}{ir}\right ] , &\\
g_2(s, r)=\left [\dfrac{f_3(s,r)}{sr}, 0, -\dfrac{f_1(s, r)}{sr}\right ], &\\
g_3(s, i)=\left [-\dfrac{f_2(s, i)}{si}, \dfrac{f_1(s, i)}{si},
0\right ], &
\end{align*}
and $f_1,\ f_2$ and $f_3$ are the right hand side of $(2-1)',\
(2-2)'$ and $(2-3)'$ reduced to functions of two variables by
using $\sum =1$ respectively. Now after some computations
\cite{Dv}, we get
$$\mbox{curl} g(s, i, r). (1,1,1)=-\left (\dfrac{b+\gamma}{is^2}+\dfrac{b}{rs^2}+
\dfrac{\beta_1}{ri^2}+\dfrac{\beta_2}{ir^2}+\dfrac{\alpha}{sr^2}\right
)< 0.$$ Since here $g.f=0$, the proof is complete by Theorem 2.1.
$\square$

\para{Corollary 2.3.} The $\omega$-limit set of any orbit of the system
$(2-1)'-(2-3)'$ with initial point in $D$ is a rest point.\\
{\bf Proof.} Since the vector field related to the system
$(2-1)'-(2-3)'$ is inward on $\partial D$ and $D$ is compact, the
$\omega$-limit set of each orbit with initial point in $D$, is a
nonempty subset of $\overset{o}{D}$. By generalized
Poincar\'e-Bendixon theorem \cite{P}, \cite{Pd} and Proposition
2.2., this set must be a rest point. $\square$

\para{Corollary 2.4.} The system $(2-1)'-(2-3)'$ has no source in
$\overset{o}{D}$.\\
{\bf Proof.}
 Suppose there is a source in $\overset{o}{D}$ for
this system. Since $D$ is positively invariant and there are
finitely many rest points in $D$ (at most four, since our system
is quadratic), there must be infinitely many heteroclinic orbit
running from this source  to another rest point. So there is a
2-gons in $\overset{o}{D}$ which is impossible by Proposition
2.2. $\square$

\para{Corollary 2.5.} Every nondegenerate rest point of the system (2-6),
(2-7) in $\overset{o}{D}_1$ is hyperbolic.\\
{\bf Proof.} Let $L$ be the linearization of this system at a rest
point $(s^{\ast}, i^{\ast})$ in $\overset{o}{D}_1$. We have to
show that if $\det L>0$, then  $\mbox{trace}\ L\neq 0$. We compute
$\mbox{trace}\ L$ at a rest point.
$$\begin{cases}
\dfrac{\partial s'}{\partial s}=(\var_2-b-\beta -\gamma )-2
\var_2 s+(\var_1-\var_2-\lambda )i,& \\[2mm]
\dfrac{\partial i'}{\partial i}=(\var_2+\var_1-\alpha -b)+
(\lambda -\var_2 )s+2(\var_1-\var_2)i.&
\end{cases}$$
From $i'=0$ and $s'=0$ at $(s^{\ast}, i^{\ast} )$, we have
$$\mbox{trace}\ L=-\dfrac{b+\gamma (1-i^{\ast} )}{s^{\ast}}-
\var_2 s^{\ast}-\dfrac{\beta_1 s^{\ast}}{i^{\ast}}+(\var_1-
\var_2)i^{\ast}.$$ If $\mbox{trace}\ L=0$, then we can slightly
increase $\var_1$ to get $\mbox{trace}\ L>0$ and determine
$\lambda$ and $\alpha$ so that $(s^{\ast}, i^{\ast})$ remains a
rest point for the new values of parameters. We may also assume
that $\det L>0$ at this rest point. Thus we obtain a source in
$\overset{o}{D}_1$ which contradicts Corollary 2.4. $\square$

\para{Remark 2.6.} A nondegenerate rest point of the system $(2-6), (2-7)$ is
obtained by a transversal intersection of two conic sections
$s'=0$ and $i'=0$. In Section 4, we shall prove that this
intersection is almost always transversal.

\section{A Special Case}

In this section we consider the planner system (2-6), (2-7) in
the case $b=\beta _1=\gamma =0$. These assumptions yields the
following system of equations

 $$\left \{ \begin{array}{lr}
s'=s(\var_2-\beta-\var_2 s+(\var_1-\var_2-\lambda)i), &
(3-1)\\
i'=i(\var_2-\var_1-\alpha+(\lambda-\var_2)s+(\var_1-\var_2)i).
 & (3-2)
\end{array}\right.$$
First of all notice that there are two invariant lines
 $s=0$ and $i=0$ with three rest points $(0, 0), \  (0, 1-
\dfrac{\alpha}{\var_2-\var_1} )$ and
$(1-\dfrac{\beta}{\var_2},0)$. The matrix of the linearized
system at $(0, 0)$ is
$$\left [
\begin{array}{cc}
\var_2-\beta & 0\\
0 & \var_2-\var_1-\alpha
\end{array}
\right ]
$$
with the eigenvalues $T_0:=\var_2-\beta$ and $T_1:=\var_2 -
\var_1 -\alpha$.

The matrix of the linearized system at $(0, 1-
\dfrac{\alpha}{\var_2 - \var_1})$ is
$$\left [
\begin{array}{cc}
\var_2-\beta+(\var_1-\var_2-\lambda )
(1-\dfrac{\alpha}{\var_2-\var_1} ) & 0\\
(\lambda-b-\var_2) (1-\dfrac{\alpha}{\var_2- \var_1}) & \var_1-
\var_2+\alpha
\end{array}
\right ]
$$
with the eigenvalues $-T_1$ and $$T_2:=\var_2-\beta +
\dfrac{(\var_1-\var_2-\lambda )(\var_2-\var_1- \alpha
)}{\var_2-\var_1}.$$ Notice that the rest point $(0,
 1-\dfrac{\alpha} {\var_2-\var_1})$ belongs to $D_1$ if
and only if $T_1\geq 0$ and coincides with $(0,0)$ in the case of
equality. It is easy to see that when $T_1\leq 0$, the origin
attracts the segment $D_1\cap\{s=0\}$ and if $T_1 >0$, then $(0,
1-\dfrac{\alpha}{\var_2 -\var_1})$  attracts
$(D_1\cap\{s=0\})-\{(0, 0)\}$.

The matrix of linearization of $(3-1),(3-2)$ at
$(1-\dfrac{\beta}{\var_2},0)$ is
\[\begin{bmatrix}
\beta-\var_2 &
(\var_1-\var_2-\lambda )(1-\dfrac{\beta}{\var_2})\\
0 & (\var_2-\var_1-\alpha)-(\var_2-\lambda)(1-\dfrac{\beta}
{\var_2})\end{bmatrix},\] with the eigenvalues $-T_0$ and
$$T_3:=\var_2-\var_1-\alpha-\dfrac{(\var_2-\lambda)(\var_2-\beta)}
{\var_2}.$$ Notice that the rest point
$(1-\dfrac{\beta}{\var_2},0)$ belongs to $D_1$ if and only if
$T_0\geq 0$ and coincides with $(0,0)$ in the case of equality.
It is easy to see that when $T_0\leq 0$, the origin attracts the
segment $D_1\cap\{i=0\}$ and if $T_1 >0$, then
$(1-\dfrac{\beta}{\var_2},0)$  attracts $(D_1\cap\{i=0\})-\{(0,
0)\}$.

\para{Proposition 3.1.} If $T_0<0$ and $T_1<0$, then the
origin is the only rest point of the system (3-1),(3-2).\\
{\bf Proof.} In this case $(0,0)$ is a sink and the rest points
$(0, 1- \dfrac{\alpha}{\var_2-\var_1} )$ and
$(1-\dfrac{\beta}{\var_2},0)$ are outside of $D_1$. Since our
planar system is quadratic, there are at most one rest point in
$\overset{o}{D}_1$. Moreover if such a rest point exists, it must
be nondegenerate, hence its Poincar\'{e} index is $\pm 1$. For a
Jordan curve $C$, let $\mu_+(C)$ and $\mu_-(C)$ denote the number
of rest points inside of $C$ with the Poincar\'{e} index +1 and
-1 respectively. Here let $C$ be the curve as shown in Figure
(3-1) which is the boundary of the union of $D_1$ and a small disk
centered at the origin. Then on this Jordan curve, our vector
field is always tangent or inward. It follows from Lemma
 5.1. (cf. Section 5) that
$I_C(X)=1$ where $X$ is the vector field related to the system
(3-1), (3-2) and $I_C(X)$ is the Poincar\'e index of $X$ with
respect to $C$. Now we can use the Poincar\'e theorem to obtain
$\mu_+(C)-\mu_-(C)=1$. Moreover we have shown that
$\mu_+(C)+\mu_-(C)\leq 2$. Therefore $\mu_+(C)=1$ and
$\mu_-(C)=0$ which finishes the proof. $\square$
\begin{center}
\unitlength 1.00mm \linethickness{0.4pt}
\begin{picture}(51.67,52.67)
\put(4.00,5.00){\line(0,1){47.33}}
\multiput(4.00,52.33)(0.12,-0.12){395}{\line(0,-1){0.12}}
\put(51.33,5.00){\line(-1,0){47.33}}
\put(4.33,1.33){\line(1,0){0.34}}
\put(22.67,24.67){\vector(-1,-1){0.20}}
\multiput(27.00,29.00)(-0.12,-0.12){37}{\line(0,-1){0.12}}
\linethickness{0.8pt} \put(4.00,21.67){\vector(0,-1){0.5}}
\put(21.00,5.00){\vector(-1,0){0.5}}
\put(8.33,5.00){\line(1,0){43.00}}
\multiput(51.33,5.00)(-0.12,0.12){395}{\line(0,1){0.12}}
\put(4.00,52.33){\line(0,-1){43.66}}
\put(4.2,4.9){\oval(8.0,8.0)[l]} \put(4.2,4.9){\oval(8.0,8.0)[rb]}
\put(4.00,4.00){\vector(0,1){0.2}}
\put(4.00,4.00){\line(0,-1){3.0}}
\put(3.00,5.00){\vector(1,0){0.2}}
\put(3.00,5.00){\line(-1,0){3.00}}
\multiput(51.67,5.00)(-0.12,0.12){398}{\line(0,1){0.12}}
\multiput(22.33,24.33)(0.12,0.12){39}{\line(0,1){0.12}}
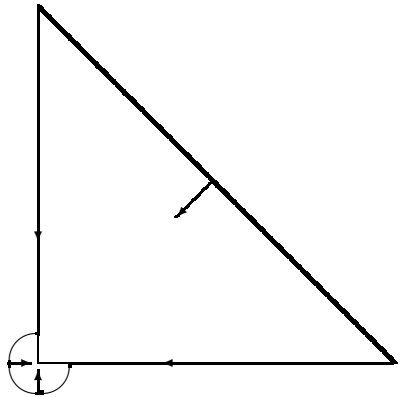
\put(30.00,-3.00){\makebox(0,0)[cc]{{\small Figure (3-1) }}}
\put(30.00,-8.00){\makebox(0,0)[cc]{{\small The Jordan curve
related to Proposition 3.1. }}}
\end{picture}
\end{center}
\vspace*{1cm}
\para{Proposition 3.2.} If $T_0>0,\ T_1>0,\ T_2<0$ and $T_3<0$, then the
system (3-1), (3-2) has a saddle point in $\overset{o}{D}_1$.\\
{\bf Proof.} These assumptions mean that the origin is a source
and $(0,1-\dfrac{\alpha}{\var_2-\var_1})$ and
$(1-\dfrac{\beta}{\var_2},0)$ are sinks. It is also concluded
that the system (3-1),(3-2) has only nondegenerate rest points.
Let $B_1, B_2$ and $B_3$ be small disks centered at the above
three rest points respectively. Let $C=\partial\Delta$ and
$\Delta =(D_1\cup B_2\cup B_3)-B_1$. (See
 Figure (3-2).) The Poincar\'e index of the Jordan curve $C$ is 1 by Lemma
 5.1. and we can use Poincar\'{e} theorem to obtain
 $\mu_+(C)-\mu_-(C)=1$.
Moreover we have  $\mu_+(C)+\mu_-(C)\leq 3$ and  $\mu_+(C)\geq 2$.
Therefore $\mu_+(C)=2$ and $\mu_-(C)=1$ which means that there is
a saddle point in $\overset{o}{D}_1$ . $\square$

\begin{center}
\unitlength=1.00mm \special{em:linewidth 0.8pt}
\special{em:ovalwidth 0.8pt} \linethickness{0.8pt}
\begin{picture}(62.34,63.33)
\put(19.67,20.00){\line(0,1){43.33}}
\put(19.67,63.33){\line(1,-1){42.67}}
\put(62.33,20.67){\line(-1,0){42.67}}
\put(19.83,46.00){\oval(7.00,6.67)[l]}
\put(20.00,21.00){\oval(7.33,7.33)[rt]}
\put(44.17,20.67){\oval(7.67,6.67)[b]}
\put(44.33,20.67){\circle*{1.49}}
\put(19.67,20.67){\circle*{1.49}}
\put(19.67,46.00){\circle*{1.49}}
\put(19.67,29.00){\vector(0,1){8.00}}
\put(19.67,60.00){\vector(0,-1){8.00}}
\put(30.00,20.67){\vector(1,0){7.00}}
\put(58.67,20.67){\vector(-1,0){7.00}}
\put(50.00,33.33){\vector(-1,-1){4.00}}
\put(22.33,23.67){\vector(1,1){4.00}}
\put(32.33,50.33){\vector(-1,-1){4.00}}
\put(44.33,16.00){\vector(0,1){3.00}}
\put(15.00,46.00){\vector(1,0){3.00}}
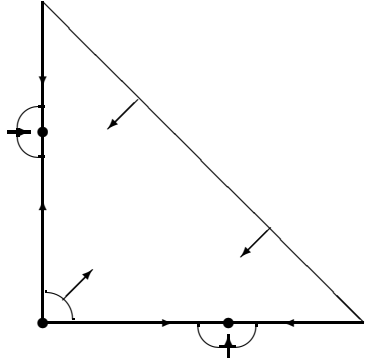
\put(40.00,9.00){\makebox(0,0)[cc]{{\small Figure (3-2)}}}
\put(40.00,3.00){\makebox(0,0)[cc]{{\small  The Jordan curve
related to Proposition 3.2.}}}
\end{picture}
\end{center}
\para{Remark 3.3.} The assumptions of the above proposition do
not contradict. To see this suppose that $T_0>0$ and $T_1>0$.
Then it is easy to check that
$$
T_2<0\Leftrightarrow
\frac{T_0}{T_1}<1+\frac{\lambda}{\var_2-\var_1}\quad ,\quad
 T_3<0\Leftrightarrow
\frac{T_0}{T_1}>(1-\frac{\lambda}{\var_2})^{-1}$$
$$(1-\frac{\lambda}{\var_2})^{-1}<(1+\frac{\lambda}{\var_2-\var_1})
\Leftrightarrow \var_1>\lambda.$$ Now if
 $\var_2>\var_1>\lambda$, one can choose $\alpha>0$
and $\beta>0$ so that $T_0>0$, $T_1>0$ and
$\dfrac{T_0}{T_1}\in((1-\dfrac{\lambda}
{\var_2})^{-1},1+\dfrac{\lambda}{\var_2-\var_1})$ to satisfy the
assumptions of Proposition 3.2. This would be helpful for the
reader who is more interested in numerical simulations.

\section{The General Case}

In this section we investigate the dynamics of the proportions
system $(2-1)'-(2-3)'$ in the general case. In order to do this,
we discuss the existence and stability of the rest points of the
planner system (2.6), (2.7). Recall that the feasibility region
$D_1$ is positively invariant and this system has no rest point on
$\partial D_1$. Indeed the vector field corresponding to the
system (2.6), (2.7) is strictly inward on $\partial D_1$. Thus by
Lemma 5.1., the Poincar\'e index of this vector field with
respect to $\partial D_1$ equals 1.  This is the first step to
prove the main result of this paper.

\para{Theorem 4.1.} If the system (2-6), (2-7) has only nondegenerate rest
points in $\overset{o}{D}_1$, then one of the following
statements holds: \\
(A) There exists a unique rest point $\overset{o}{D}_1$ which is a sink and
attracts $D_1(GAS)$. \\
(B) There are two sinks and a saddle point in
$\overset{o}{D}_1$.\\
Moreover both of them occur for suitable values of the involved
parameters. \\
{\bf Proof.} All rest points in $\overset{o}{D}_1$ are hyperbolic
by Corollary 2.5. Let $\mu_0$, $\mu_1$ and $\mu_2$ be the number
of sinks, saddles and sources in $\overset{o}{D}_1$. Since the
index of $\partial D_1$ is 1, we have $\mu_0-\mu_1+\mu_2=1$ in
$D_1$. Furthermore $\mu_0+\mu_1+\mu_2\leq 4$ and $\mu_2=0$ by
Corollary 2.4. Thus we have either $\mu_0=1, \mu_1=\mu_2=0$ or
$\mu_0=2, \mu_1=1,
 \mu_2=0$. Now the first conclusion gives (A) and the second one gives
(B). In order to see that each of the two above statements occurs,
recall that in the case of Theorem 3.1., there is neither a saddle
nor a nonhyperbolic rest point in $D_1$. Thus for small values of
$b_1$ and $\gamma$, there cannot be any saddle point or
nonhyperbolic rest point in $\overset{o}{D}_1$. This means that
$(A)$ occurs. Similarly, under small perturbation, the saddle
point obtained in Theorem 3.2 remains in $\overset{o}{D}_1$ and
yields the case $(B)$ of this theorem.  $\square$

\para{Remark 4.2.} It is well-known that the basin of attraction
 of a sink is a connected open set. In the case (B) of the
above result,  $D_1$ is the distinct union of the basins of
attraction of these two sinks and the stable manifold of the
saddle point in $D_1$. The basins of attraction are open and
$D_1$ is connected. Thus the stable manifold of the saddle point
separates them. In order to specify these sets, we can
numerically find those two points at which $\omega$-limit set
changes from one sink to another when someone moves on $\partial
D_1$. (See Figure 4.1.)

\begin{center}
\unitlength=1.00mm \special{em:linewidth 0.8pt}
\linethickness{0.8pt}
\begin{picture}(66.00,72.00)
\put(19.00,25.00){\line(0,1){47.00}}
\put(19.00,72.00){\line(1,-1){47.00}}
\put(66.00,25.00){\line(-1,0){47.00}}
\put(26.00,51.67){\circle*{1.89}}
\put(31.00,40.00){\circle*{1.89}}
\put(43.00,32.67){\circle*{2.00}}
\put(20.00,55.00){\vector(3,-2){4.00}}
\put(25.33,61.00){\vector(0,-1){6.00}}
\put(33.67,54.33){\vector(-3,-1){6.00}}
\put(45.67,40.67){\vector(-1,-3){2.00}}
\put(54.33,28.67){\vector(-3,1){8.33}}
\put(27.00,31.00){\vector(1,0){9.33}}
\bezier{200}(39.00,47.67)(33.00,38.33)(28.00,48.67)
\put(30.00,45.33){\vector(-2,3){2.33}}
\put(38.00,26.00){\vector(1,1){4.33}}
\bezier{200}(20.00,38.00)(24.00,37.67)(23.67,48.00)
\put(23.67,43.67){\vector(0,1){4.33}}
\put(54.00,35.33){\vector(-4,-1){7.67}}
\bezier{300}(41.67,42.67)(35.67,37.67)(40.67,35.00)
\bezier{300}(21.67,34.33)(28.67,37.67)(33.33,33.00)
\put(39.33,36.00){\vector(1,-1){1.00}}
\put(31.67,34.33){\vector(3,-2){2.00}}
\bezier{104}(41.00,48.00)(38.67,41.33)(20.67,36.33)
\bezier{68}(37.00,35.00)(32.33,37.00)(26.00,47.33)
\put(37.00,42.67){\vector(-3,-1){2.00}}
\put(26.67,38.00){\vector(2,1){2.00}}
\put(27.00,45.67){\vector(-2,3){0.67}}
\put(35.33,35.67){\vector(3,-1){1.67}}
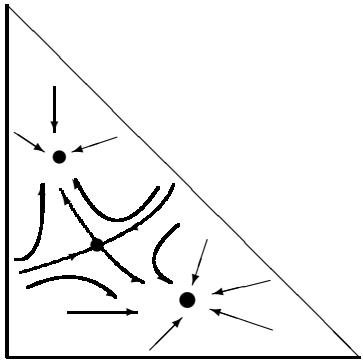
\put(41.00,16.00){\makebox(0,0)[cc]{{\small Figure (4-1)}}}
\put(41.00,10.00){\makebox(0,0)[cc]{{\small The phase portrait of
the case (B) of Theorem 4.1.}}}
\end{picture}
\end{center}

\para{Remark 4.3.} The non-uniqueness of endemic equilibrium
proportions yields some interesting conclusions. The most
significant one is that the initial condition of the population
may also be important besides the involved parameters. The effect
of the initial condition is more crucial when the population
equations (2-1)-(2-4) is considered. From (2-4) we get
$$\frac{N'}{N}=b-d-\var_1 i-\var_2 r.$$ Now suppose that a
solution $(s(t), i(t),r(t))$ of the system $(2-1)'-(2-3)'$ tends
to an equilibrium $(s^*,i^*,s^*)$. If we set
$T=\dfrac{b}{d+\var_1i^*+\var_2 r^*}$, then $N(t)\rightarrow
\infty$ if $T>1$ and $N(t)\rightarrow 0$ if $T<1$. (See \cite{Bv1}
for more details.) Now each endemic equilibrium gives a $T$ and
when there are two endemic equilibria, we may get different
values for $T$ at these two points.

It remains to consider the case in which there is a degenerate
rest point in $D_1$ for our planar system. Let $\Omega$ be the
parameters space of the system (2-6), (2-7) as an open subset of
$\R^8_+$ and $\Omega_1$ be the set of all possible values of
parameters for which the system (2-6), (2-7) has a nonhyperbolic
(or equivalently degenerate by Corollary 2.5.) rest point in
$D_1$. The following fact about $\Omega_1$ shows that our problem
has fairly been solved.

\para{Proposition 4.4.} With the above notations, $\Omega_1$ is a closed
nonempty subset of $\Omega$ with zero measure.\\
{\bf Proof.} We first show that $\Omega_1$ is closed and nonempty.
Since $D_1$ is compact and all rest points in the statements (A)
and (B) of Theorem 4.1 are hyperbolic, both (A) and (B) occur in
open subsets of $\Omega$. (In other words, our system is
structurally stable in the nondegenerate case.) Since $\Omega$ is
connected, it cannot be the union of these two distinct open
subsets. Thus $\Omega_1$ is closed, nonempty and indeed large
enough to separate two open subsets. We use Sard's theorem
\cite{GP} to show that $\Omega_1$ has zero measure. Notice that
from (2-6), we can write $i$ in terms of $s$ if $\gamma
+(\var_2-\var_1+\lambda )s\neq 0$ and from (2-7) we can write $s$
in terms of $i$ if $\beta_1+(\lambda -\var_2)i\neq 0$. If
$\beta_1+(\lambda -\var_2)i=0$ then from $i'=0$ in (2-7), we have
$(\var_2-\var_1-b-\alpha )+ (\var_1-\var_2)i=0$. Thus $\var_2-
\var_1-b-\alpha >0$ and hence $\gamma +(\var_2-\var_1+ \lambda
)s>0$. It follows that either $i$ can be written in terms of $s$
from (2-6) or $s$ in terms of $i$ from (2-7). In the first case
we have a root for the equation $h_1(s)=\beta_1$ and in the
latter case a root for $h_2(i)=\gamma$ where: $$\begin{cases}
h_1(s)=(b+\var_1-\var_2+\alpha)\dfrac{i}{s}+(\var_2
-\lambda)i+(\var_2-\var_1)\dfrac{i^2}{s},&\\
h_2(i)=-\dfrac{b}{r}+(b+\beta)\dfrac{s}{r}+ (\lambda
-\var_1)\dfrac{is}{r}-\var_2s,&
\end{cases}$$
and $r=1-s-i$ as before. Notice that in the first equation $i$ has
be written in terms of $s$ and in the second equation $s$ has been
written in terms of $i$. Since at a degenerate rest point of the
system (2-6), (2-7), the curves $i'=0$ and $s'=0$ are not
transverse, it makes $\beta_1$ a critical value of $h_1$ or
$\gamma$ a critical value of $h_2$. In other words, whenever the
system (2-6),(2-7) has a degenerate rest point in $D_1$, either
$\beta_1$ or $\lambda$ belongs to the set of critical values
which has zero measure by Sard's theorem. Therefore  $\Omega_1$
is contained in the union of two sets with zero measure. $\square$

\section{The Index Lemma}

A basic fact which has been used during the proof of our results
is that the Poincar\'e index of a piece-wise smooth Jordan curve
on which the vector field is either tangent or inward is always 1.
Here we provide the proof of this fact. The reader is referred to
 \cite{P} for more details about the Poincar\'e index.

Let $\gamma :[0,1]\rg\R^2$
 be a piecewise smooth Jordan curve i.e. there exists
a sequence $0=t_0<t_1<\ldots<t_n=1$ such that $\gamma$
 is smooth on $(t_i,t_{i+1})$
for $0\leq i\leq n-1$. Suppose $\gamma '$ is always nonzero on
$(t_i,t_{i+1})$, moreover the left and right derivatives of
$\gamma$ at $t_i$ exist and both  are nonzero. With these
assumptions we can define the external angle $\theta_i\in (-\pi
,\pi )$ at $t_i$. Also the inward normal vector $N(t)$ is defined
for $t\not=t_i$ and its right and left limit $N(t_i^{\pm})$ exist
at each $t_i$. For such a curve $\gamma$ we prove the following
lemma.

\para{Lemma 5.1.} Let $U$ be a neighborhood of the image of $\gamma$ and
$X:U\rg \R^2$ be a smooth vector field which does not vanish on
the image of $\gamma$ and satisfies $N(t)\cdot
X(\gamma(t))\geq 0$ for $t\not=t_i$. Then $I_{\gamma}(X)=1$.\\
{\bf Proof.} We define a new curve $\gamma_1:[0,n+1]\rg\R^2$ by
\begin{align*}
&\gamma_1(t)=\left\{\begin{array}{lll} \gamma(t-i) & i+t_i\leq
t\leq i+t_{i+1}, & (0\leq i\leq n-1)\\ \gamma(t_{i+1}) &
i+t_{i+1}<t<i+1+t_{i+1}. &
\end{array}\right.
\end{align*}
Notice that $\gamma_1$
 moves like $\gamma$ but stops at $\gamma(t_i)$ for a unit of
 time, hence $I_{\gamma}(X)=I_{\gamma_1}(X)$.
 Now we use this unit of time to rotate $N(t_i^-)$
to arrive at $N(t^+_i)$. To do this, we define a continuous
function $N_1:[0,n+1]\rg\R^2$ by
\begin{align*}
&N_1(t)=\left\{\begin{array}{ll} N(t-i) & i+t_i< t< i+t_{i+1},\\
R_{t\theta_i}(N(t_i^-)) & i+t_{i+1}\leq t\leq i+1+t_{i+1},
\end{array}\right.
\end{align*}
where $R_{\alpha}$ is the rotation function with the angle
$\alpha$. Now both $X(\gamma_1(t))$ and $N_1(t)$ are continuous
and do not vanish on $[0,n+1]$.

\noindent{\bf Claim:} $N_1(t)\cdot X(\gamma_1(t))\geq 0$.\\
Since $ N(t)\cdot X(\gamma(t))\geq 0\ \mbox{for}\ t\neq t_i$, by
continuity we have $N(t_i^{\pm})\cdot X(\gamma (t_i))\geq 0$.
Moreover $\theta_i\in (-\pi ,\pi )$, thus for any $t\in [0, 1]$,
there exists  $t'\in [0,1]$ such that
$$R_{t\theta_i}(N(t^-_i))=\frac{t'N(t^-_i)+(1-t')N(t^+_i)}{||t'N(t^-_i)+
(1-t')N(t^+_i)||}.$$ So $R_{t\theta_i}(N(t^-_i)). X(\gamma
(t_i))\geq 0$. As a result of the above claim
$$
|\Delta\Theta (X(\gamma_1(t)))-\Delta\Theta N_1(t)|\leq\pi
 \eqno{(5-1)}
$$
where $\Delta\Theta$ is the total change of the angle on $[0,
n+1]$. Then we can write
\begin{align*}
\Delta\Theta (X(\gamma_1(t)))&=\Delta\Theta (X(\gamma (t)))=2\pi I_{\gamma}(X).\\
\Delta\Theta (N_1(t))&=\sum\limits_{i=0}^{n-1}\Delta\Theta
(N_1|_{(i+t_i, i+t_{i+1})})+\sum\limits_{i=0}^{n-1}\Delta\Theta
(N_1|_{(i+t_{i+1}, i+1+t_{i+1}
)})\\
&=\sum\limits_{i=0}^{n-1}\Delta\Theta (\gamma '|_{(i+t_i,
i+t_{i+1})})+\sum\limits_{i=1}^{n}\theta i\\
&=2\pi .
\end{align*}
by Gauss-Bonnet theorem. Now from (5-1) we get $|2\pi
I_{\gamma}(X)-2\pi |<\pi$, hence $|I_{\gamma}(X)-1|<
\dfrac{1}{2}.$ Since $I_{\gamma}(X)$ is an integer, we get
$I_{\gamma}(X)=1.\ \square$

\para{\bf Acknowledgement.} The author would like to thank
Institute for Studies in Theoretical Physics and Mathematics for
supporting this research.

\end{document}